\pdfoutput=1
\RequirePackage{ifpdf}
\ifpdf % We are running pdfTeX in pdf mode
\documentclass[pdftex]{sigma}
\else
\documentclass{sigma}
\fi

\begin{document}

\allowdisplaybreaks

\renewcommand{\PaperNumber}{046}

\FirstPageHeading

\ShortArticleName{On Addition Formulae for Sigma Functions of Telescopic Curves}

\ArticleName{On Addition Formulae for Sigma Functions\\
of Telescopic Curves}

\Author{Takanori AYANO~$^\dag$ and Atsushi NAKAYASHIKI~$^\ddag$}

\AuthorNameForHeading{T.~Ayano and A.~Nakayashiki}

\Address{$^\dag$~Department of Mathematics, Osaka University, Toyonaka, Osaka 560-0043, Japan}
\EmailD{\href{mailto:tayano7150@gmail.com}{tayano7150@gmail.com}}

\Address{$^\ddag$~Department of Mathematics, Tsuda College, Kodaira, Tokyo 187-8577, Japan}
\EmailD{\href{mailto:atsushi@tsuda.ac.jp}{atsushi@tsuda.ac.jp}}

\ArticleDates{Received March 13, 2013, in f\/inal form June 14, 2013; Published online June 19, 2013}

\Abstract{A telescopic curve
is a certain algebraic curve def\/ined by $m-1$ equations in the af\/f\/ine
space of dimension $m$, which can be a hyperelliptic curve and
an $(n,s)$ curve as a~special case.
We extend the  addition formulae for sigma functions of $(n,s)$ curves to those of telescopic curves. The expression of the prime form in terms
of the derivative of the sigma function is also given.}


\Classification{14H70; 37K20; 14H55; 14K25}

\section{Introduction}

In this paper we study the multivariate sigma functions of telescopic curves and derive addition formulae together with their degenerate limits.

\looseness=1
The multivariate sigma function originally introduced by F.~Klein \cite{Kl1,Kl2} for hyperelliptic curves is generalized and
extensively studied for the last decade (see~\cite{BEL3} and references therein). Compared with Riemann's theta function,
the sigma function is more algebraic and is directly related with
the def\/ining equation of an algebraic curve. A typical example where this
nature of the sigma function is exhibited is the inversion problem
of algebraic integrals. It is well known that the solution of Jacobi's inversion problem for a hyperelliptic curve has a simple description
 by hyperelliptic $\wp$-functions, the second logarithmic
derivatives of the sigma function~\cite{Bak,BEL3}. The inversion of
 hyperelliptic or more general algebraic integrals of genus $g$ on the Abel--Jacobi image~$W_k$ of the $k$-th symmetric products of the curve
with $k<g$ is extensively studied in connection with the problem of mathematical physics (see \cite{BEF, BEL3} and references therein). This problem is intimately related with the problem on the vanishing of the derivatives
of the sigma function on $W_k$.
Recently it is recognized that the approach from the view point
of tau functions of integrable hierarchies provides a general and
ef\/fective method to study such a prob\-lem~\cite{NY}.

In Sato's theory of Kadomtsev--Petviashvili (KP) hierarchy \cite{SS} the tau function
is constructed from a point of the universal Grassmann manifold (UGM).
For a solution corresponding to an algebraic curve the point of UGM is
specif\/ied by expanding functions or sections of bundles on the curve
using a local
coordinate at a given point. The point of UGM obtained in this way belongs
to the cell of UGM labeled by the partition $\lambda$ determined from the
gap sequence at the point. The series expansion of the corresponding tau function begins from Schur function associated with~$\lambda$.
The tau function corresponding to the point of UGM specif\/ied by the af\/f\/ine
ring of an $(n,s)$ curve~\cite{BEL2} had been used to study the sigma function in~\cite{N2,NY}.

Therefore it is important to have such a pair
$(X,p_\infty)$ of an algebraic curve $X$ and a point $p_\infty\in X$
that satisf\/ies the following two conditions.
The f\/irst is that a basis, as a vector space, of the space of regular functions on $X\backslash \{p_\infty\}$ can be explicitly described.  The second is that the gap sequence at $p_\infty$ can be computable.
A traditional example is a non-singular plane algebraic curve which can be completed by one point $\infty$ $(=p_\infty)$ such as a hyperelliptic curve of odd degree or more generally an $(n,s)$ curve.
Telescopic curves give new examples.  They can be hyperelliptic and $(n,s)$ curves as special cases and are not realized as non-singular plane algebraic curves in general.
Before explaining telescopic curves let us brief\/ly explain the origin of the term ``telescopic''.

Let $a_1,\dots,a_m$ be relatively prime positive integers.
For a nonnegative integer $n$ the problem of determining nonnegative integers $x_1,\dots,x_m$ satisfying the Diophantine equation
\begin{gather}
n=a_1x_1+\cdots+a_mx_m\label{dio}
\end{gather}
has been studied in number theory since early times.
It is well known that the equation~\eqref{dio} has a solution if $n$ is suf\/f\/iciently large.
The greatest number $n$ for which the equation~\eqref{dio} has no solution is called Frobenius number and we denote it by $F(a_1,\dots,a_m)$.
Brauer \cite{B1} gave an upper bound of the Frobenius number as
\begin{gather}
F(a_1,\dots,a_m)\le -a_1+\sum_{i=2}^ma_i\left(\frac{d_{i-1}}{d_i}-1\right),\label{tele}
\end{gather}
where $d_i={\rm gcd}(a_1,\dots ,a_i)$.
Also Brauer \cite{B1} and Brauer and Seelbinder \cite{B2} showed that the equality in~\eqref{tele} holds if and only if
\begin{gather}
\frac{a_i}{d_i}\in \frac{a_1}{d_{i-1}}{\mathbb Z}_{\geq 0}+
\cdots+\frac{a_{i-1}}{d_{i-1}}{\mathbb Z}_{\geq 0},
\qquad
2\leq i\leq m.\label{telesco}
\end{gather}
The condition~\eqref{telesco} is introduced in Brauer \cite{B1} for the f\/irst time and now it is called telescopic condition.
For $a_1,\dots,a_m$ satisfying~\eqref{telesco} the semigroup $S:=a_1\mathbb{Z}_{\ge0}+\cdots+a_m\mathbb{Z}_{\ge0}$ is called telescopic semigroup.
The telescopic semigroup has many nice structures and has many applications in algebraic geometric code, algebraic curve cryptography, and commutative algebra (see for instance \cite{K, M, Miu}).

In \cite{Miu} Miura introduced a certain canonical form, Miura canonical form, for def\/ining equations of any non-singular algebraic curve.
A telescopic curve \cite{Miu} is a special curve for which Miura canonical form is easy to determine.
Let $m\geq 2$ and $(a_1,\dots ,a_m)$ a sequence of relatively prime
positive integers satisfying the telescopic condition~\eqref{telesco}.
Then the telescopic curve associated with $(a_1,\dots ,a_m)$ or the $(a_1,\dots ,a_m)$ curve is the algebraic curve def\/ined by
certain $m-1$ equations in ${\mathbb C}^m$.
The form of def\/ining equations is explicitly computable from $(a_1,\dots ,a_m)$ (see~\eqref{eq-2-2}).
If a telescopic curve is non-singular, then it can be completed by
adding one point, say $\infty$, and the gap sequence at~$\infty$
becomes the complement of the telescopic semigroup~\cite{Aya1,Miu}.
In such a case the genus of the curve is also explicitly computable (see~\eqref{eq-2-4}).

In \cite{NY} the vanishing and the expansion of the sigma function
of an $(n,s)$ curve on the Abel--Jacobi image $W_k$ for $k<g$
are studied using the properties of the tau function of the
KP-hierarchy. Those results are then applied to study the restriction of the addition formulae on $W_k$ to the lower strata $W_{k'}$ with $k'<k$.
On the other hand the sigma function of the telescopic curve is explicitly constructed in~\cite{Aya1}.
In this paper we show that almost results in~\cite{NY} are extended to the
case of telescopic curves. The results imply two things. The f\/irst is
that telescopic curves are natural objects to study sigma functions.
We expect, more generally, the Miura canonical form is suitable to
describe properties of the sigma functions. The second is that the
tau function approach is ef\/fective in a more general case than that of
$(n,s)$ curves.
 We expect that the method by integrable hierarchies can equally be ef\/f\/icient to study sigma functions with characteristics of arbitrary Riemann surfaces~\cite{KS}.

Finally we comment that the sigma functions of certain space curves,
which are not telescopic, are studied in \cite{Komeda,Ma}.

The present paper is organized as follows.
In Section~\ref{section2} the def\/inition and examples of telescopic curves are given.
The construction of the sigma function, up to the normalization constant, associated with telescopic curves is reviewed in Section~\ref{section3}. The local coordinate $z$ at $\infty$ is specif\/ied and the expression in terms of~$z$ of the variables appearing in the def\/ining equations of the curve is given.
This is necessary to determine constants appearing in every formula
 in later sections.
In Section~\ref{section4} the expression of the tau function is given using the sigma function. The normalization constant necessary in the def\/inition of the sigma function is specif\/ied with the help of it.
 In Section~\ref{section5} main results including the addition formulae are given.
Their proofs are indicated in Section~\ref{section6}.
 The example of addition formulae is given for a $(4,6,5)$
curve in Section~\ref{section7}.
In Appendix~\ref{appendixA} the detailed properties on the series expansion of the sigma function
are given for the sake of completeness of the construction of
the sigma function.

\section{Telescopic curves}\label{section2}

In this section we brief\/ly review the def\/inition and properties of
telescopic curves following \cite{Aya1, Miu} and give some examples.

Let $m\geq 2$, $(a_1,\dots ,a_m)$ a sequence of positive integers
such that ${\rm gcd}(a_1,\dots ,a_m)=1$ and $d_i={\rm gcd}(a_1,\dots ,a_i)$ for
$1\leq i\leq m$. We call $(a_1,\dots ,a_m)$ telescopic if
\[
\frac{a_i}{d_i}\in \frac{a_1}{d_{i-1}}{\mathbb Z}_{\geq 0}+
\cdots+\frac{a_{i-1}}{d_{i-1}}{\mathbb Z}_{\geq 0},
\qquad
2\leq i\leq m.
\]

The following examples of telescopic sequences are given in \cite{Nishi}.

\begin{example}\label{example1}\qquad
\begin{enumerate}\itemsep=0pt
\item[$(i)$] $(a_1,a_2)$, s.t.\ ${\rm gcd}(a_1,a_2)=1$.

\item[$(ii)$] $(k,k+2,k+1)$, s.t.\ $k$ even.

\item[$(iii)$] $(ab,bc,a+c)$, s.t.\ ${\rm gcd}(a,c)=1$, ${\rm gcd}(b,a+c)=1$.

\item[$(iv)$] $(a_1,\dots ,a_m)$, s.t.\ $a_i=a^{m-i}b^{i-1}$, $a>b$, ${\rm gcd}(a,b)=1$.
\end{enumerate}
\end{example}

Notice that whether a sequence is telescopic depends on the order
of the numbers. For example, $(4,6,5)$ is telescopic while $(4,5,6)$ is not.

In the following we assume that $A_m=(a_1,\dots ,a_m)$ is telescopic
unless otherwise stated.

Let
\[
B(A_m)=\left\{(l_1,\dots ,l_m)\in {\mathbb Z}_{\geq 0}^m\,|\,0\leq l_i\leq \frac{d_{i-1}}{d_i}-1\  \text{for}\ 2\leq i\leq m\right\}.
\]

\begin{lemma}[\cite{Aya1, Miu}]\label{x}
For any $a\in a_1\mathbb{Z}_{\ge0}+\cdots+a_m\mathbb{Z}_{\ge0}$, there exists a unique element $(k_1,\dots,k_m)$ of $B(A_m)$ such that
\[
\sum_{i=1}^ma_ik_i=a.
\]
\end{lemma}

By this lemma, for any $2\leq i\leq m$, there exists a unique sequence
$(l_{i1},\dots ,l_{im})\in B(A_m)$ satisfying
\begin{gather}
\sum_{j=1}^m a_j l_{ij}= a_i\frac{d_{i-1}}{d_i}.
\label{eq-2-3}
\end{gather}
Consider $m-1$ polynomials in $m$ variables $x_1,\dots ,x_m$ given by
\begin{gather}
F_i(x)=x_i^{d_{i-1}/d_i}-\prod_{j=1}^m x_j^{l_{ij}}-
\sum \kappa^{(i)}_{j_1\dots j_m}x_1^{j_1}\cdots x_m^{j_m},
\quad
2\leq i\leq m,
\label{eq-2-2}
\end{gather}
where $\kappa^{(i)}_{j_1\dots j_m}\in\mathbb{C}$
and the sum of the right hand side is over all
$(j_1,\dots ,j_m)\in B(A_m)$ such that
\[
\sum_{k=1}^m a_kj_k<a_i\frac{d_{i-1}}{d_i}.
\]

Let $X^{\rm{af\/f}}$ be the common zeros of $F_2$,\dots ,$F_ m$:
\[
X^{\rm{af\/f}}=\big\{(x_1,\dots ,x_m)\in\mathbb{C}^m\,|\,F_i(x_1,\dots ,x_m)=0,\, 2\leq i\leq m\big\}.
\]
In \cite{Aya1, Miu} $X^{\rm{af\/f}}$ is proved to be an af\/f\/ine algebraic curve.
We assume that $X^{\rm{af\/f}}$ is nonsingular. Let $X$ be the compact Riemann surface corresponding to $X^{\rm{af\/f}}$. Then $X$ is obtained from $X^{\rm{af\/f}}$
by adding one point, say $\infty$ \cite{Aya1,Miu}.
It is proved in \cite{Aya1,Miu} that $x_i$ has a pole of order $a_i$ at
$\infty$.
The genus of $X$ is given by \cite{Aya1,Miu}
\begin{gather}
g=\frac{1}{2}\left(
1+\sum_{i=2}^m a_i\frac{d_{i-1}}{d_i}-\sum_{i=1}^m a_i
\right).
\label{eq-2-4}
\end{gather}
We call $X$ the $(a_1,\dots ,a_m)$ curve or the telescopic curve
associated with $(a_1,\dots ,a_m)$.
The numbers $a_1$,\dots ,$a_m$ are a generator of the semigroup of
non-gaps at $\infty$.

\begin{example}\label{example2}\qquad
\begin{enumerate}\itemsep=0pt
\item[$(i)$] The telescopic curve associated with a pair of relatively prime
integers $(n,s)$ is the $(n,s)$ curve introduced in \cite{BEL2}.

\item[$(ii)$] For $A_3=(2k,2k+2,2k+1)$, $k\geq 2$, in $(ii)$ of Example~\ref{example1},
polynomials $F_i$ are given by
\begin{gather*}
F_2(x)=x_2^k-x_1^{k+1}-
\sum\nolimits^{(2)}\kappa^{(2)}_{i_1,i_2,i_3}x_1^{i_1}x_2^{i_2}x_3^{i_3},
\\
F_3(x)=x_3^2-x_1x_2-\sum\nolimits^{(3)}\kappa^{(3)}_{i_1,i_2,i_3}x_1^{i_1}x_2^{i_2}x_3^{i_3},
\end{gather*}
where $\sum^{(i)}$, $i=2,3$ signify the sum
over all $(i_1,i_2,i_3)\in B(A_3)$ such that
\[
 2ki_1+2(k+1)i_2+(2k+1)i_3<
\begin{cases}
2k(k+1)&\   \text{for $\sum\nolimits^{(2)}$},\\
2(2k+1)&\   \text{for $\sum\nolimits^{(3)}$}.
\end{cases}
\]
The genus of $X$ is $g=k^2$.

\item[$(iii)$] For $A_3=(ab,bc,a+c)$, $a\neq 1$, in $(iii)$ of Example~\ref{example1}, we have
\begin{gather*}
F_2(x)=x_2^a-x_1^c-\sum\nolimits^{(2)}\kappa^{(2)}_{i_1,i_2,i_3}x_1^{i_1}x_2^{i_2}x_3^{i_3},
\\
F_3(x)=x_3^b-x_1x_2-\sum\nolimits^{(3)}\kappa^{(3)}_{i_1,i_2,i_3}x_1^{i_1}x_2^{i_2}x_3^{i_3},
\end{gather*}
where $\sum^{(i)}$, $i=2,3$ denote the sum over all
$(i_1,i_2,i_3)\in B(A_3)$ such that
\[
abi_1+bci_2+(a+c)i_3<
\begin{cases}
abc&\ \text{for $\sum\nolimits^{(2)}$},\\
b(a+c)&\ \text{for $\sum\nolimits^{(3)}$}.
\end{cases}
\]
The genus of $X$ is
\[
g=\frac{1+abc-a-c}{2}.
\]

\item[$(iv)$] For $A_m=(a_1,\dots ,a_m)$ in $(iv)$ of Example~\ref{example1}, we have
\[
F_i(x)=x_i^a-x_{i-1}^b-
\sum_{a_1j_1+\cdots+a_mj_m<aa_i}
\kappa^{(i)}_{j_1\dots j_m}x_1^{j_1}\cdots x_m^{j_m}.
\]
The genus of $X$ is
\[
g=\frac{a-b+(b-1)a^m-(a-1)b^m}{2(a-b)}.
\]
\end{enumerate}
\end{example}

\section{Sigma function of telescopic curves}\label{section3}

An algebraic bilinear dif\/ferential of a telescopic curve associated with
$(a_1,\dots ,a_m)$ is explicitly constructed in \cite{Aya1}. Consequently
an expression of the sigma function in terms of Riemann's theta function
and some algebraic data had been given.
We recall the results of \cite{Aya1} and add some necessary
results for our purpose.

Let $X$ be a telescopic curve of genus $g\geq 1$
associated with $(a_1,\dots ,a_m)$ and $(l_{i1},\dots ,l_{im})$ the element
of $B(A_m)$ specif\/ied by~\eqref{eq-2-3}.

\begin{lemma}\label{a}
For any $i$ we have $l_{ij}=0$ for $j\geq i$.
\end{lemma}

\begin{proof}
Since $A_m$ is telescopic, there exist $k_1,\dots,k_{i-1}\in\mathbb{Z}_{\ge0}$ such that
\begin{gather}
a_i\frac{d_{i-1}}{d_i}=a_1k_1+\cdots+a_{i-1}k_{i-1}.
\label{lemma2-eq1}
\end{gather}
We prove that we can take $0\leq k_j< d_{j-1}/d_j$ for any $j\ge2$ by changing
$k_j$ appropriately if necessary.

Suppose that $k_{j'}\ge d_{j'-1}/d_{j'}$ for some $j'$. Take the largest number $j$ satisfying this condition.
Let us write
\[
k_j=\frac{d_{j-1}}{d_j}q+r,
\]
with $q\ge1$, $0\le r<d_{j-1}/d_j$.
Since $A_m$ is telescopic, there exist $u_1,\dots, u_{j-1}\in\mathbb{Z}_{\ge0}$ such that
\[
a_j\frac{d_{j-1}}{d_j}=a_1u_1+\cdots+a_{j-1}u_{j-1}.
\]
Then we have
\[
a_jk_j=a_j\frac{d_{j-1}}{d_j}q+a_jr=a_1qu_1+\cdots+a_{j-1}qu_{j-1}+a_jr.
\]
Substituting this into~\eqref{lemma2-eq1} we get the expression
of the form~\eqref{lemma2-eq1} with $0\leq k_{l}<d_{l-1}/d_l$ for
any $l\geq j$. Repeating similar change of $k_{j'}$ for $j'$ smaller than $j$ successively we f\/inally get the expression of $a_id_{i-1}/d_i$
of the form~\eqref{lemma2-eq1} with $0\leq k_j<d_{j-1}/d_j$
for any $j\geq 2$.

By the def\/inition of $B(A_m)$, $(k_1,\dots,k_{i-1},0,\dots,0)\in B(A_m)$.
Since the element of $B(A_m)$ satisfying~\eqref{eq-2-3} is unique by Lemma~\ref{x}, $(l_{i1},\dots ,l_{im})=(k_1,\dots,k_{i-1},0,\dots,0).$
\end{proof}

For the def\/ining equations~\eqref{eq-2-2},
we assign degrees as
\[
\deg \kappa_{j_1\dots j_m}^{(i)}=a_id_{i-1}/d_i-\sum_{k=1}^ma_kj_k.
\]

\begin{lemma}\label{b}
It is possible to take a local parameter $z$ around $\infty$
such that
\begin{gather}
x_1=\frac{1}{z^{a_1}},
\qquad
x_k=\frac{1}{z^{a_k}}\left(1+\sum_{l=1}^{\infty}e_{kl}z^l\right),\qquad 2\leq k\leq m,
\label{eq-3-1}
\end{gather}
where $e_{kl}$ belongs to $\mathbb{Q}\big[\big\{\kappa^{(i)}_{j_1\dots j_m}\big\}\big]$ and is homogeneous of degree $l$ if $e_{kl}\neq0$.
\end{lemma}

\begin{proof}
It is possible to take a local parameter $z_0$ around $\infty$ such that
\[
x_1=\frac{1}{z_0^{a_1}}.
\]
Let $\zeta=\exp\big(2\pi\sqrt{-1}/a_1\big)$ and $i\ge0$.
Then $z_i:=\zeta^iz_0$ is also a local parameter around $\infty$.
Let~$e_k^{(i)}$ be the coef\/f\/icient of the f\/irst term of the series expansion of $x_k$ around $\infty$ with respect to $z_i$:
\begin{gather}
x_k=\frac{e_k^{(i)}}{z_i^{a_k}}(1+O(z_i)),\qquad 2\le k\le m.
\label{h}
\end{gather}
We prove that there exists $i$ such that $e_2^{(i)}=\cdots=e_m^{(i)}=1$.

Let $e^{(i)}=\big(e_2^{(i)},\dots,e_m^{(i)}\big)$ for $0\le i<a_1$.
First we show $e^{(i)}\neq e^{(j)}$ for $i\neq j$.
Suppose $e^{(i)}=e^{(j)}$.
Since $e_k^{(i)}=\zeta^{a_ki}e_k^{(0)}$, we have $\zeta^{a_k(i-j)}=1$ for $k=2,\dots,m$.
From ${\rm gcd}(a_1,\dots,a_m)=1$ and $0\le i,j<a_1$, we have $i=j$.

By Lemma~\ref{a} the def\/ining equations of $X$ are as follows:
\begin{gather}
x_k^{d_{k-1}/d_k}=x_1^{l_{k1}}\cdots x_{k-1}^{l_{k k-1}}+\sum \kappa^{(k)}_{j_1\dots j_m}x_1^{j_1}\cdots x_m^{j_m},
\qquad
2\leq k\leq m.
\label{i}
\end{gather}
By substituting~\eqref{h} to~\eqref{i} and comparing the coef\/f\/icients of $z_i^{-a_kd_{k-1}/d_k}$, we have
\[
\big(e_2^{(i)}\big)^{d_1/d_2}=1,
\qquad
\big(e_k^{(i)}\big)^{d_{k-1}/d_k}=\big(e_2^{(i)}\big)^{l_{k2}}\cdots \big(e_{k-1}^{(i)}\big)^{l_{k k-1}},
\qquad
3\le k\le m.
\]

Let
\[
S=\big\{(s_2,\dots,s_m)\in\mathbb{C}^{m-1}\;|\;s_2^{d_1/d_2}=1, \;s_k^{d_{k-1}/d_k}=s_2^{l_{k2}}\cdots s_{k-1}^{l_{k k-1}},\;3\le k\le m\big\}.
\]
Since $\sharp S=(d_1/d_2)\cdots(d_{m-1}/d_m)=(d_1/d_m)=a_1$ and $e^{(i)}\in S$ for $i=0,\dots,a_1-1$, we have
\[
S=\big\{e^{(0)},\dots,e^{(a_1-1)}\big\}.
\]
Since $(1,\dots,1)\in S$, there exists $i$ such that $e^{(i)}=(1,\dots,1)$.
For $z:=z_i$, $x_k$ is expanded as
\[
x_1=\frac{1}{z^{a_1}},\qquad x_k=\frac{1}{z^{a_k}}\left(1+\sum_{l=1}^{\infty}e_{kl}z^l\right),\qquad e_{kl}\in\mathbb{C}.
\]

Let us prove that $e_{kl}$ belongs to $\mathbb{Q}\big[\big\{\kappa^{(i)}_{j_1\dots j_m}\big\}\big]$ and is homogeneous of degree $l$ if $e_{kl}\neq0$.
We def\/ine the order $<$ in the set $\{e_{kl}\}$ so that $e_{k'l'}<e_{kl}$ if
\begin{enumerate}\itemsep=0pt
\item[1)] $l'<l$ or
\item[2)] $l'=l$ and $k'<k$.
\end{enumerate}

We prove the statement by induction on this order.
By~\eqref{eq-2-2} and Lemma~\ref{a} we have
\begin{gather}
\left(1+\sum_{j=1}^{\infty}e_{kj}z^j\right)^{\frac{d_{k-1}}{d_k}}
=\prod_{s=2}^{k-1}
\left(1+\sum_{j=1}^{\infty}e_{sj}z^j\right)^{l_{ks}}\nonumber\\
\hphantom{\left(1+\sum_{j=1}^{\infty}e_{kj}z^j\right)^{\frac{d_{k-1}}{d_k}} =}{}
+\sum\kappa_{j_1\dots j_m}^{(k)}z^{\frac{a_kd_{k-1}}{d_k}-\sum\limits_{s=1}^ma_sj_s}
\prod_{s=2}^m\left(1+\sum\limits_{j=1}^{\infty}e_{sj}z^j\right)^{j_s},
\label{gty}
\end{gather}
where we def\/ine the empty product from $s=2$ to $1$ to be one
in the f\/irst term of the right hand side.

In~\eqref{gty} for $k=2$, the coef\/f\/icient of $z$ of the left hand side is $(d_1/d_2)e_{21}$ and that of the right hand side is the sum of $\kappa_{j_1\dots j_m}^{(2)}$ with $(j_1,\dots,j_m)$ satisfying the
equation $(a_2d_1/d_2)-\sum\limits_{s=1}^ma_sj_s=1$.
Therefore the statement is valid for the minimal element~$e_{21}$.

Assume that the statement holds for any $e_{k'l'}$ satisfying $e_{k'l'}<e_{kl}$.
The coef\/f\/icient of $z^l$ of the left hand side of~\eqref{gty} is
$(d_{k-1}/d_k)e_{kl}+T$, where $T$ is a sum of $\prod_i e_{kq_i}$ satisfying $\sum_iq_i=l$ and $q_i<l$.
In the right hand side of~\eqref{gty},
the coef\/f\/icient of $z^l$ of the f\/irst term is the sum of $\prod_ie_{p_iq_i}$ satisfying $2\le p_i<k$ and $\sum_iq_i=l$,
and that of the second term
is the sum of $\kappa_{j_1\dots j_m}^{(k)}\prod_ie_{p_iq_i}$ with $(j_1,\dots ,j_m)$ satisfying $\sum_iq_i=l-(a_kd_{k-1}/d_k)+\sum\limits_{s=1}^ma_sj_s$.
Therefore, by the assumption of induction, we f\/ind that $e_{kl}$ belongs to $\mathbb{Q}\big[\big\{\kappa^{(i)}_{j_1\dots j_m}\big\}\big]$ and is homogeneous of degree $l$ if $e_{kl}\neq0$.
\end{proof}

For a meromorphic function $f$ on $X$ we denote by ${\rm ord}_{\infty}(f)$
the order of a pole at $\infty$. Then we have ${\rm ord}_{\infty}(x_i)=a_i$.
We enumerate the monomials $x_1^{\alpha_1}\cdots x_m^{\alpha_m},\;(\alpha_1,\dots,\alpha_m)\in B(A_m)$ according as the order of a pole at $\infty$ and denote them by
$\varphi_i$, $i\geq 1$. In particular we have $\varphi_1=1$.

Let $(w_1,\dots ,w_g)$ be the gap sequence at $\infty$:
\[
\{w_i\,|\,1\leq i\leq g\}={\mathbb Z}_{\geq0}
\Big\backslash\left\{\sum_{i=1}^m a_i{\mathbb Z}_{\geq0}\right\},
\qquad
w_1<\cdots<w_g.
\]
In particular $w_1=1$, since $g\geq 1$.

A basis of holomorphic one forms is given by
\begin{gather}
du_{w_i}=-\frac{\varphi_{g+1-i}}{\det G(x)}dx_1,\label{eq-3-2}
\end{gather}
where $G(x)$ is the Jacobian matrix
\[
G(x)=\left(\frac{\partial F_i}{\partial x_j}\right)_{2\leq i,j\leq m}.
\]
The following lemma is proved in \cite{Aya1}.

\begin{lemma}\label{lem-3-1}
We have $w_{g}=2g-1$.
 In particular $du_{2g-1}$ has a zero of order $2g-2$ at $\infty$.
\end{lemma}

More precisely we have the following properties.

\begin{proposition}\label{c}\qquad
\begin{enumerate}\itemsep=0pt
\item[$(i)$] The following expansion is valid around $\infty$:
\[
du_{2g-1}=z^{2g-2}\left(1+\sum_{l=1}^\infty e_l'z^l\right)dz,
\]
where $e_l'$ belongs to ${\mathbb Q}\big[\big\{\kappa^{(i)}_{j_1\dots j_m}\big\}\big]$ and
is homogeneous of degree $l$ if $e_l'\neq 0$.

\item[$(ii)$] For $1\leq i\leq g$ the expansion of $du_{w_i}$ at $\infty$ is of the form
\[
du_{w_i}=z^{w_i-1}(1+O(z))dz.
\]
\end{enumerate}
\end{proposition}

\begin{proof}
$(i)$
From Lemmas~\ref{a}, \ref{b}, we have, around $\infty$,
\[\det G(x)=a_1z^{-\sum\limits_{i=2}^m((d_{i-1}/d_i)-1)a_i}\left(1+\sum_{l=1}^{\infty}e_l''z^l\right)dz,\]
where $e_l''$ belongs to $\mathbb{Q}[\{\kappa^{(i)}_{j_1\dots j_m}\}]$ and is homogeneous of degree $l$ if $e_l''\neq0$.
Therefore, from~\eqref{eq-2-4}, we obtain the assertion.

$(ii)$
Let $w_i^*={\rm ord}_{\infty}(\varphi_i)$, $W=\{w_1,\dots,w_g\}$, and $W'=\{w_1^*,\dots,w_g^*\}$.
Note that $W\cup W'=\{0,1,\dots,2g-1\}$.
If $w\in W'$, then $w_g-w\in W$.
In fact, if $w_g-w\in W'$, then $w_g\in W'$, which is contradiction.
Since $w_g=2g-1$, we have $2g-1-w_{g+1-i}^*=w_i$ for any $i$.
Therefore, from~\eqref{eq-3-2}, Lemma~\ref{b}, and Proposition~\ref{c}$(i)$, we obtain the assertion.
\end{proof}

The algebraic bilinear dif\/ferential takes the form
\[
\widehat{\omega}(x,y)=d_y\Omega(x,y)+
\sum_{i=1}^g du_{w_i}(x)dr_i(y),
\]
where $x=(x_1,\dots ,x_m)$, $y=(y_1,\dots ,y_m)$ are points on $X$,
\[
\Omega(x,y)=\frac{\det H(x,y)}{(x_1-y_1)\det G(x)}dx_1,
\]
$H=(h_{ij})_{2\leq i,j\leq m}$ with
\[
h_{ij}=\frac{F_i(y_1,\dots ,y_{j-1},x_j,x_{j+1},\dots ,x_m)-F_i(y_1,\dots ,y_{j-1},y_j,x_{j+1},\dots ,x_m)}{x_j-y_j},
\]
and $dr_i$ is a second kind dif\/ferential with a pole only at
$\infty$. By construction $\{du_{w_i},dr_i\}$ becomes a symplectic
basis of the cohomology group $H^1(X,{\mathbb C})$ (see \cite{Aya1, N1}).

Take a symplectic basis $\{\alpha_i,\beta_i\}$ of the homology group
and def\/ine the period matrices by
\begin{gather*}
2\omega_1=\left(\int_{\alpha_j}du_{w_i}\right),
\qquad
2\omega_2=\left(\int_{\beta_j}du_{w_i}\right),
\\
-2\eta_1=\left(\int_{\alpha_j}dr_i\right),
\qquad
-2\eta_2=\left(\int_{\beta_j}dr_i\right).
\end{gather*}
The normalized period matrix is given by $\tau=\omega_1^{-1}\omega_2$.

Let $\delta=\tau\delta'+\delta''$, $\delta',\delta''\in {\mathbb R}^g$
be the Riemann's constant with respect to the choice
$(\{\alpha_i,\beta_i\},\infty)$. We set $\delta={}^t({}^t\delta',
{}^t\delta'')$.
Since $du_{2g-1}$ has a zero of order $2g-2$ at $\infty$ by Lemma \ref{lem-3-1}, we have $\delta\in(\mathbb{Z}/2)^{2g}$.

We def\/ine the function $\widehat{\sigma}(u)$, $u=(u_{w_1},\dots ,u_{w_g})$
by
\[
\widehat{\sigma}(u)=\exp\left(\frac{1}{2}{}^tu\eta_1\omega_1^{-1}u\right)
\theta[\delta]\big((2\omega_1)^{-1}u,\tau\big),
\]
where $\theta[\delta](u)$ is the Riemann's theta function with the characteristic $\delta$.

\section{Relation with tau function}\label{section4}

In the case of $(n,s)$ curves the expression of the tau function
of the KP-hierarchy in terms of the sigma function is given in \cite{EEG,EH, N2}. For the tau function corresponding to the point
of Sato's universal Grassmann manifold (UGM) specif\/ied by the af\/f\/ine ring
of a telescopic curve a similar expression holds.

Let $A$ be the af\/f\/ine ring of $X$,
\[
A={\mathbb C}[x_1,\dots ,x_m]/I,
\]
where $I$ is the ideal generated by $F_i(x)$, $2\leq i\leq m$.
As a vector space $A=\oplus_{i=1}^\infty {\mathbb C} \varphi_i$.

We embed $A$ in UGM as in Section 5 of~\cite{N2} using the local parameter $z$ in~\eqref{eq-3-1} and denote~$U^A$ the
image of it. Let $\xi^A$ be the normalized frame of $U^A$ and
$\tau(t;\xi^A)$ the tau function corresponding to $\xi^A$. Then
we have the expansion of the form
\begin{gather}
\tau(t;\xi^A)=s_\lambda(t)+\sum_{\lambda<\mu}\xi_\mu s_\mu(t),
\label{eq-4-1}
\end{gather}
where $\lambda=(\lambda_1,\dots ,\lambda_g)$ is the partition def\/ined
by
\begin{gather}
\lambda_i=w_{g+1-i}-g+i,
\qquad
1\leq i\leq g,
\label{eq-4-1-1}
\end{gather}
$s_\mu(t)$ is the Schur function corresponding to the partition $\mu$
and, in general, for two partitions $\mu=(\mu_1,\dots ,\mu_l)$ and
$\nu=(\nu_1,\dots ,\nu_l)$, $\mu\leq \nu$ if and only if
$\mu_i\leq \nu_i$ for any $i$.

Def\/ine $b_{ij}$, $\widehat{q}_{ij}$ and $c_i$ by the expansions
\begin{gather*}
du_{w_i}=\sum_{j=1}^\infty b_{ij}z^{j-1}dz,
\qquad
\widehat{\omega}(p_1,p_2)=\left(
\frac{1}{(z_1-z_2)^2}+\sum_{i,j=1}^\infty \widehat{q}_{ij}z_1^{i-1}z_2^{j-1}\right)dz_1dz_2,
\\
\log\left( z^{-g+1}\sqrt{\frac{du_{2g-1}}{dz}}\right)=
\sum_{i=1}^\infty c_i\frac{z^i}{i},
\end{gather*}
where $z_i=z(p_i)$.\footnote{In the right hand side of the
def\/ining equation of $c_i$ in~\cite{N2}, $z_i$ should be corrected
as $z^i/i$.}
Let us set
\[
B=(b_{ij})_{1\leq i\leq g, 1\leq j},
\qquad
\widehat{q}(t)=\sum_{i,j=1}^\infty \widehat{q}_{ij}t_it_j,
\qquad
t={}^t(t_1,t_2,\dots ).
\]
Then we have

\begin{theorem}\label{th-4-1}\qquad
\begin{enumerate}\itemsep=0pt
\item[$(i)$] There exists a constant $C$ such that
\[
\tau(t;\xi^A)=
C\exp\left(-\sum_{i=1}^\infty c_it_i+\frac{1}{2}\widehat{q}(t)\right)
\widehat{\sigma}(Bt).
\]

\item[$(ii)$] The tau function $\tau(t;\xi^A)$ is a solution to the $a_1$-reduced
KP-hierarchy.
\end{enumerate}
\end{theorem}

\begin{proof}
The proof of this theorem is completely pararell to
the case of $(n,s)$ curves in~\cite{N2}. In fact the special property
of $(n,s)$ curves we use in~\cite{N2}
is the existence of a holomorphic
one form which vanishes at $\infty$ to the order $2g-2$. In the present
case $du_{2g-1}$ has such a property by Lemma~\ref{lem-3-1}.
\end{proof}

Combining the expansion~\eqref{eq-4-1} with Theorem~\ref{th-4-1} we have

\begin{corollary}\label{cor-4-1}
We have the following expansion
\[
C\widehat{\sigma}(u)=s_{\lambda}(t)|_{t_{w_i}=u_{w_i}}+\cdots,
\]
where $\lambda$ is defined by~\eqref{eq-4-1-1} and $\cdots$ part
is a series in $\prod\limits_{i=1}^g {u_{w_i}}^{\gamma_i}$,
$\sum\limits_{i=1}^g \gamma_i w_i>\sum\limits_{i=1}^g \lambda_i$.
\end{corollary}

\begin{definition} We def\/ine the sigma function by
\[
\sigma(u)=C\widehat{\sigma}(u).
\]
\end{definition}

It follows from this def\/inition and Corollary~\ref{cor-4-1} that the series
expansion of $\sigma(u)$ at the origin begins from Schur function
corresponding to the gap sequence at~$\infty$.
It is possible to give more precise properties of the expansion of~$\sigma(u)$ which is similar to the case of $(n,s)$ curves (see Appendix~\ref{appendixA}).

\section{Addition formulae}\label{section5}
Our result is that all properties of the sigma functions of $(n,s)$ curves
given in Sections~4 and~5 of the paper~\cite{NY} are valid, formally without any change, for sigma functions of telescopic curves. The strategy of the
proofs of theorems in this section is explained in the next section.

In order to state the results precisely
we need the prime function of a telescopic curve which was introduced
in \cite{N1} for $(n,s)$ curves.

Let $E(p_1,p_2)$ be the prime form of a telescopic curve $X$.
Since $du_{2g-1}$ has a zero of order $2g-2$ at $\infty$
by Lemma~\ref{lem-3-1}, it
is possible to def\/ine, as in the case of $(n,s)$ curves, the prime function ${\tilde E}(p_1,p_2)$  by
\[
{\tilde E}(p_1,p_2)=-E(p_1,p_2)\prod_{i=1}^2\sqrt{du_{2g-1}(p_i)}
\exp\left(\frac{1}{2}\int_{p_1}^{p_2}{}^td{\mathbf u}\eta_1\omega_1^{-1}
\int_{p_1}^{p_2}d{\mathbf u}\right).
\]
It is a multi-valued analytic function on $X\times X$ which has similar  properties to that for $(n,s)$ curves.

For a partition $\lambda=(\lambda_1,\dots ,\lambda_l)$ and $0\leq k\leq l$
we set
\begin{gather*}
(w_l,\dots ,w_1)=(\lambda_1+l-1,\lambda_2+l-2,\dots ,\lambda_l),
\\
N_{\lambda,k}=\lambda_{k+1}+\cdots+\lambda_l,\qquad
N'_{\lambda,1}=\lambda_2+\cdots+\lambda_l-l+1,
\\
c'_{\lambda,k}=\frac{N_{\lambda,k}!}{\prod\limits_{i=1}^{l-k}w_i!}
\prod_{i<j}^{l-k}(w_j-w_i),
\qquad
\tilde{c}_{\lambda}=\frac{N'_{\lambda,1}!}{\prod\limits_{i=1}^{l-1}(w_i-1)!}
\prod_{i<j}^{l-1}(w_j-w_i),
\end{gather*}
where $c'_{\lambda,l}$ is considered to be $1$.

The following theorems give the properties of the sigma function
restricted to the Abel--Jacobi image $W_k$ of the $k$-th symmetric products of $X$ for $k<g$.

\begin{theorem}\label{th-5-1}
Let $1\leq k\leq g$ and $p_1,\dots ,p_k\in X$. Then
\begin{enumerate}\itemsep=0pt
\item[$(i)$] We have
\[
\partial_{u_1}^{N_{\lambda,k}}\sigma\left(\sum_{i=1}^k p_i\right)=
c'_{\lambda,k}S_{(\lambda_1,\dots ,\lambda_k)}(z_1,\dots ,z_k)+\cdots,
\]
where $\cdots$ part is a series in $z_1,\dots ,z_k$ containing only terms
proportional to $\prod\limits_{i=1}^k z_i^{\gamma_i}$ with
$\sum\limits_{i=1}^k \gamma_i>\sum\limits_{i=1}^k\lambda_i$.

\item[$(ii)$] The following expansion in $z_k$ holds:
\[
\partial_{u_1}^{N_{\lambda,k}}
\sigma\left( \sum_{i=1}^k p_i\right)
=
\frac{c'_{\lambda,k}}{c'_{\lambda,k-1}}
\partial_{u_1}^{N_{\lambda,k-1}}\sigma
\left(\sum_{i=1}^{k-1} p_i\right)
z_k^{\lambda_k}
+O\big(z_k^{\lambda_k+1}\big).
\]
\end{enumerate}
\end{theorem}

\begin{theorem}\label{th-5-2}\qquad
\begin{enumerate}\itemsep=0pt
\item[$(i)$] If $n<N'_{\lambda,1}$ we have, for $p_1,p_2\in X$,
\[
\partial_{u_1}^n\sigma(p_1-p_2)=0.
\]

\item[$(ii)$] The following expansion with respect to $z_i=z(p_i)$, $i=1,2$
is valid:
\[
\partial_{u_1}^{N'_{\lambda,1}} \sigma(p_1-p_2)
=(-1)^{g-1} \tilde{c}_{\lambda} (z_1z_2)^{g-1}(z_1-z_2)(1+\cdots),
\]
where $\cdots$ part is a series in $z_1$, $z_2$ which
contains only terms proportional to $z_1^iz_2^j$ with $i+j>0$.

\item[$(iii)$] We have
\[
\partial_{u_1}^{N'_{\lambda,1}}\sigma(p_1-p_2)=
(-1)^{g-1}
\frac{\tilde{c}_\lambda}{c'_{\lambda,1}}
\partial_{u_1}^{N_{\lambda,1}}\sigma(p_1)z_2^{g-1}+O\big(z_2^g\big).
\]
\end{enumerate}
\end{theorem}

The next theorem gives the expression of the prime function as a derivative
of the sigma function.

\begin{theorem}\label{th-5-3}
Let $\lambda=(\lambda_1,\dots ,\lambda_g)$ be the partition corresponding
to the gap sequence at $\infty$ of an $(a_1,\dots ,a_m)$ curve $X$.
Then
\[
{\tilde E}(p_1,p_2)=(-1)^{g-1}\tilde{c}_\lambda^{-1}\partial_{u_1}^{N'_{\lambda,1}}
\sigma(p_1-p_2).
\]
\end{theorem}

\begin{corollary}\label{cor-5-1}
For $p\in X$ we have
\[
{\tilde E}(\infty,p)=
{c'_{\lambda,1}}^{-1}\partial_{u_1}^{N_{\lambda,1}}\sigma(p).
\]
\end{corollary}

Finally the addition formulae of sigma functions for telescopic
curves are given by

\begin{theorem}\label{the-4}\qquad
\begin{enumerate}\itemsep=0pt
\item[$(i)$] For $n\geq g$ and $p_i\in X$, $1\leq i\leq n$,
\[
\frac{\sigma\left( \sum\limits_{i=1}^n p_i\right)
\prod\limits_{i<j}\partial_{u_1}^{N'_{\lambda,1}}
\sigma(p_j-p_i)}{\prod\limits_{i=1}^n\big(\partial_{u_1}^{N_{\lambda,1}}\sigma(p_i)\big)^n}
=\tilde{b}_{\lambda,n}
\det\left(\varphi_{i}(p_j)\right)_{1\leq i,j\leq n},
\]
with
\[
\tilde{b}_{\lambda,n}=
(-1)^{\frac{1}{2}gn(n-1)}\tilde{c}_{\lambda}^{\frac{1}{2}n(n-1)}(c'_{\lambda,1})^{-n^2}.
\]

\item[$(ii)$] For $n<g$
\[
\frac{\partial_{u_1}^{N_{\lambda,n}}
\sigma\left( \sum\limits_{i=1}^n p_i\right)
\prod\limits_{i<j}\partial_{u_1}^{N'_{\lambda,1}}\sigma(p_j-p_i)}
{\prod\limits_{i=1}^n\big(\partial_{u_1}^{N_{\lambda,1}}\sigma(p_i)\big)^n}
=b_{\lambda,n}'
\det\left(\varphi_{i}(p_j)\right)_{1\leq i,j\leq n},
\]
with
\[
b_{\lambda,n}'=
(-1)^{\frac{1}{2}gn(n-1)}\tilde{c}_{\lambda}^{\frac{1}{2}n(n-1)}(c'_{\lambda,1})^{-n^2}
c'_{\lambda,n}.
\]
\end{enumerate}
\end{theorem}

The addition formulae of this kind were f\/irstly derived by \^{O}nishi~\cite{O} for the hyperelliptic sigma functions.
We remark that the formulae in this theorem are written
 using algebraic data only, the sigma function, its derivatives
 and the algebraic functions~$\varphi_i$.
This fact makes it possible to study the restriction of the
formulae on~$W_k$ to lower strata~$W_{k'}$ for~$k'<k$ as in $(ii)$.
This is a main dif\/ference
of our formulae and Fay's general addition formulae~\cite{F}.

\section{Proofs}\label{section6}

In \cite{NY} properties of the sigma function
of an $(n,s)$ curve have been proved by establishing the corresponding
properties of Schur and tau functions. In this paper we exclusively
consider
$t_1$-derivatives and $u_1$-derivatives. We omit the results on
``$a$-derivatives'' in~\cite{NY}, since they are not used in addition
formulae for general telescopic curves.

As far as $t_1$-derivatives
 are concerned, all the statements for Schur and tau functions
in~\cite{NY} hold, as stated there, for Schur function~$s_\lambda(t)$ associated with any partition~$\lambda$ and the tau functions~$\tau(t)$ which have the expansion of the form
\[
\tau(t)=s_\lambda(t)+\sum_{\lambda<\mu}\xi_\mu s_\mu(t).
\]

Then Theorems~\ref{th-5-1}, \ref{th-5-2}, \ref{th-5-3}, Corollary~\ref{cor-5-1} can be proved in a similar manner to the case of~$(n,s)$
curves using Theorem~\ref{th-4-1}.

By Theorem~\ref{th-5-3} the prime function of a telescopic curve can be written as a derivative of the sigma function.
Conversely the sigma function can be expressed by using the prime function
and algebraic functions $\varphi_i$ as in the case of an
$(n,s)$ curve.

\begin{theorem}\label{th-6-1}
For $n\geq g$ and $p_i\in X$, $1\leq i\leq n$,
\begin{gather}
\sigma\left( \sum_{i=1}^n p_i\right)=
\frac{\prod\limits_{i=1}^n {\tilde E}(\infty,p_i)^n}{\prod\limits_{i<j}^n{\tilde E}(p_i,p_j)}
\det\left(\varphi_i(p_j)\right)_{1\leq i,j\leq n}.
\label{eq-th61}
\end{gather}
\end{theorem}

Expanding~\eqref{eq-th61} in $z(p_n)$ successively with the
help of Theorem~\ref{th-5-1} we get

\begin{corollary}\label{cor-6-1}
For $n<g$ we have
\[
\partial_{u_1}^{N_{\lambda,n}}
\sigma\left( \sum_{i=1}^n p_i\right)=
c'_{\lambda,n}
\frac{\prod\limits_{i=1}^n {\tilde E}(\infty,p_i)^n}{\prod\limits_{i<j}^n{\tilde E}(p_i,p_j)}
\det\left(\varphi_i(p_j)\right)_{1\leq i,j\leq n}.
\]
\end{corollary}

Theorem~\ref{the-4} can be obtained from Theorem~\ref{th-6-1} and Corollary~\ref{cor-6-1}
by substituting the sigma function expression of the prime function
given by Theorem~\ref{th-5-3} and Corollary~\ref{cor-5-1}.

\section[Example: $(4,6,5)$-curve]{Example: $\boldsymbol{(4,6,5)}$-curve}\label{section7}

In this section we give an explicit example of the addition formulae in the case of a $(4,6,5)$-curve $X$.
By Example~\ref{example2}$(ii)$ in Section~\ref{section2}, the genus of $X$ is~4.
The gap sequence at $\infty$ is $(w_1,w_2,w_3,w_4)=(1,2,3,7)$.
The partition corresponding to the gap sequence at $\infty$ is $\lambda=(\lambda_1,\lambda_2,\lambda_3,\lambda_4)=(4,1,1,1)$.
Therefore we have $N_{\lambda,0}=7$, $N_{\lambda,1}=3$, $N_{\lambda,2}=2$, $N_{\lambda,3}=1$, $N_{\lambda,1}'=0$, $c_{\lambda,1}'=c_{\lambda,2}'=c_{\lambda,3}'=\tilde{c}_{\lambda}=1$.
On the other hand we have $\varphi_1=1$, $\varphi_2=x_1$, $\varphi_3=x_3$, $\varphi_4=x_2$.
Therefore the addition formulae given in Theorem~\ref{the-4} for $n=2,3,4$ are as follows:

\begin{enumerate}\itemsep=0pt
\item[$(i)$] For $n=2,$
\[
\frac{\left(\partial_{u_1}^2\sigma(p_1+p_2)\right)
\sigma(p_2-p_1)}{(\partial_{u_1}^3\sigma(p_1))^2(\partial_{u_1}^3\sigma(p_2))^2}
=x_1(p_2)-x_1(p_1).
\]

\item[$(ii)$] For $n=3,$
\[
\frac{\partial_{u_1}\sigma(p_1+p_2+p_3)
\prod\limits_{1\leq i<j\leq 3}\sigma(p_j-p_i)}
{\prod\limits_{j=1}^3(\partial_{u_1}^3\sigma(p_j))^3}
=\begin{vmatrix} 1 & 1 & 1  \\ x_1(p_1) & x_1(p_2) & x_1(p_3) \\ x_3(p_1) & x_3(p_2) & x_3(p_3) \end{vmatrix}.
\]

\item[$(iii)$] For $n=4,$
\[
\frac{\sigma(p_1+p_2+p_3+p_4)
\prod\limits_{1\leq i<j\leq 4}\sigma(p_j-p_i)}
{\prod\limits_{j=1}^4(\partial_{u_1}^3\sigma(p_j))^4}
=\begin{vmatrix} 1 & 1 & 1 & 1 \\ x_1(p_1) & x_1(p_2) & x_1(p_3) & x_1(p_4) \\ x_3(p_1) & x_3(p_2) & x_3(p_3) & x_3(p_4) \\ x_2(p_1) & x_2(p_2) & x_2(p_3) & x_2(p_4)\end{vmatrix}.
\]
\end{enumerate}

\appendix

\section{Series expansion of sigma function}\label{appendixA}

Using Theorem~\ref{th-6-1} in a similar manner to the case of $(n,s)$ curves (cf.~\cite{N1}), we can show the following theorem.

\begin{theorem}\label{main}
The expansion of $\sigma(u)$ at the origin takes the form
\[
\sigma(u)=s_{\lambda}(t)|_{t_{w_i}=u_{w_i}}+\sum_{\sum\limits_{i=1}^g\gamma_iw_i>
\sum\limits_{i=1}^g\lambda_i}\tilde{e}_{\gamma_1\dots\gamma_g}u_{w_1}^{\gamma_1}\cdots u_{w_g}^{\gamma_g},
\]
where $\tilde{e}_{\gamma_1\dots\gamma_g}$ belongs to $\mathbb{Q}\big[\big\{\kappa^{(i)}_{j_1\dots j_m}\big\}\big]$ and is homogeneous of degree ${\sum\limits_{i=1}^g\!\gamma_iw_i-\!\sum\limits_{i=1}^g\!\lambda_i}$ if $\tilde{e}_{\gamma_1\dots\gamma_g}\!\neq\!0$.
\end{theorem}

\begin{remark}
It is possible to prove the above theorem using the relation with the tau function as in \cite{N2}.
\end{remark}

\subsection*{Acknowledgements}

The authors would like to thank the referees for the useful comments.
This research was partially supported  by Grant-in-Aid for JSPS Fellows (22-2421) and for Scientif\/ic Research (C) 23540245 from Japan Society for the Promotion of Science.

\pdfbookmark[1]{References}{ref}
\LastPageEnding


\begin{thebibliography}{99}
\footnotesize\itemsep=0pt

\bibitem{Aya1}
Ayano T., Sigma functions for telescopic curves, \textit{Osaka~J. Math.}, {t}o
  appear, \href{http://arxiv.org/abs/1201.0644}{arXiv:1201.0644}.

\bibitem{Bak}
Baker H.F., Abelian functions. Abel's theorem and the allied theory of theta
  functions, \textit{Cambridge Mathematical Library}, Cambridge University Press,
  Cambridge, 1995.

\bibitem{BEF}
Braden H.W., Enolski V.Z., Fedorov Yu.N., Dynamics on strata of trigonal
  Jacobians and some integrable problems of rigid body motion,
  \href{http://arxiv.org/abs/1210.3596}{arXiv:1210.3596}.

\bibitem{B1}
Brauer A., On a problem of partitions, \textit{Amer.~J. Math.} \textbf{64}
  (1942), 299--312.

\bibitem{B2}
Brauer A., Seelbinder B.M., On a problem of partitions.~{II}, \textit{Amer.~J.
  Math.} \textbf{76} (1954), 343--346.

\bibitem{BEL3}
Buchstaber V.M., Enolski V.Z., Leykin D.V., Multi-dimensional sigma functions,
  \href{http://arxiv.org/abs/1208.0990}{arXiv:1208.0990}.

\bibitem{BEL2}
Bukhshtaber V.M., Enolski V.Z., Leykin D.V., Rational analogues of abelian
  functions, \href{http://dx.doi.org/10.1007/BF02465189}{\textit{Funct. Anal. Appl.}} \textbf{33} (1999), 83--94.

\bibitem{EEG}
Eilbeck J.C., Enolski V.Z., Gibbons J., Sigma, tau and {A}belian functions of
  algebraic curves, \href{http://dx.doi.org/10.1088/1751-8113/43/45/455216}{\textit{J.~Phys.~A: Math. Theor.}} \textbf{43} (2010),
  455216, 20~pages, \href{http://arxiv.org/abs/1006.5219}{arXiv:1006.5219}.

\bibitem{F}
Fay J.D., Theta functions on {R}iemann surfaces, \textit{Lecture Notes in
  Mathematics}, Vol.~352, Springer-Verlag, Berlin, 1973.

\bibitem{EH}
Harnad J., Enolski V.Z., Schur function expansions of KP $\tau$-functions
  associated to algebraic curves, \href{http://dx.doi.org/10.1070/RM2011v066n04ABEH004755}{\textit{Russ. Math. Surv.}} \textbf{66}
  (2011), 767--807, \href{http://arxiv.org/abs/1012.3152}{arXiv:1012.3152}.

\bibitem{K}
Kirfel C., Pellikaan R., The minimum distance of codes in an array coming from
  telescopic semigroups, \href{http://dx.doi.org/10.1109/18.476245}{\textit{IEEE Trans. Inform. Theory}} \textbf{41}
  (1995), 1720--1732.

\bibitem{Kl1}
Klein F., Ueber hyperelliptische {S}igmafunctionen, \href{http://dx.doi.org/10.1007/BF01445285}{\textit{Math. Ann.}}
  \textbf{27} (1886), 431--464.

\bibitem{Kl2}
Klein F., Ueber hyperelliptische {S}igmafunctionen (Zweite Abhandlung),
  \href{http://dx.doi.org/10.1007/BF01443606}{\textit{Math. Ann.}} \textbf{32} (1888), 351--380.

\bibitem{Komeda}
Komeda J., Matsutani S., Previato E., The sigma function for Weierstrass
  semigroups $\langle 3,7,8 \rangle$ and $\langle 6,13,14,15,16 \rangle$,
  \href{http://arxiv.org/abs/1303.0451}{arXiv:1303.0451}.

\bibitem{KS}
Korotkin D., Shramchenko V., On higher genus {W}eierstrass sigma-function,
  \href{http://dx.doi.org/10.1016/j.physd.2012.01.002}{\textit{Phys.~D}} \textbf{241} (2012), 2086--2094, \href{http://arxiv.org/abs/1201.3961}{arXiv:1201.3961}.

\bibitem{Ma}
Matsutani S., Sigma functions for a space curve (3,4,5) type with an appendix
  by J.~Komeda, \href{http://arxiv.org/abs/1112.4137}{arXiv:1112.4137}.

\bibitem{M}
Micale V., Olteanu A., On the {B}etti numbers of some semigroup rings,
  \href{http://dx.doi.org/10.4418/2012.67.1.13}{\textit{Matematiche (Catania)}} \textbf{67} (2012), 145--159,
  \href{http://arxiv.org/abs/1111.1433}{arXiv:1111.1433}.

\bibitem{Miu}
Miura S., Linear codes on af\/f\/ine algebraic curves, \textit{Trans. IEICE}
  \textbf{J81-A} (1998), 1398--1421.

\bibitem{N1}
Nakayashiki A., On algebraic expressions of sigma functions for {$(n,s)$}
  curves, \textit{Asian~J. Math.} \textbf{14} (2010), 175--211,
  \href{http://arxiv.org/abs/0803.2083}{arXiv:0803.2083}.

\bibitem{N2}
Nakayashiki A., Sigma function as a tau function, \href{http://dx.doi.org/10.1093/imrn/rnp135}{\textit{Int. Math. Res. Not.}}
  \textbf{2010} (2010), no.~3, 373--394, \href{http://arxiv.org/abs/0904.0846}{arXiv:0904.0846}.

\bibitem{NY}
Nakayashiki A., Yori K., Derivatives of Schur, tau and sigma functions on
  Abel--Jacobi images, in Symmetries, Integrable Systems and Representations,
  \href{http://dx.doi.org/10.1007/978-1-4471-4863-0_17}{\textit{Springer Proceedings in Mathematics \& Statistics}}, Vol.~40, Editors
  K.~Iohara, S.~Morier-Genoud, B.~Remy, Springer-Verlag, London, 2013,
  429--462, \href{http://arxiv.org/abs/1205.6897}{arXiv:1205.6897}.

\bibitem{Nishi}
Nishijima D., Order counting algorithm for Fermat curves and Klein curves using
  $p$-adic cohomology, Master's Thesis, Osaka University, 2008.

\bibitem{O}
{\^O}nishi Y., Determinant expressions for hyperelliptic functions (with an
  Appendix by Shigeki Matsutani), \href{http://dx.doi.org/10.1017/S0013091503000695}{\textit{Proc. Edinb. Math. Soc.~(2)}}
  \textbf{48} (2005), 705--742, \href{http://arxiv.org/abs/math.NT/0105189}{math.NT/0105189}.

\bibitem{SS}
Sato M., Sato Y., Soliton equations as dynamical systems on
  inf\/inite-dimensional {G}rassmann manifold, in Nonlinear Partial Dif\/ferential
  Equations in Applied Science ({T}okyo, 1982), \textit{North-Holland Math.
  Stud.}, Vol.~81, Editors P.D.~Lax, H.~Fujita, G.~Strang, North-Holland,
  Amsterdam, 1983, 259--271.

\end{thebibliography}
\end{document}